\documentclass{article}
\usepackage[utf8]{inputenc}
\usepackage[english]{babel}
\usepackage{amsthm,amsmath,amsbsy,amsfonts}
\usepackage{color,graphicx}

\usepackage[colorlinks=true,urlcolor=blue,
citecolor=red,linkcolor=blue,linktocpage,pdfpagelabels,
bookmarksnumbered,bookmarksopen]{hyperref} 

\title{Semilinear elliptic problems in $\mathbb{R}^N$: the interplay between the potential and the nonlinear term}
\author{Elves Alves de Barros e Silva  \\
	Departamento de Matem\'atica, \\Universidade de Bras\'{\i}lia \\
Bras\'{i}lia,  70910-900, Brazil \\
	\and
	Sergio H. Monari Soares \\
	Instituto de Ci\^encias Matem\'aticas e de Computa\c{c}\~ao\\
         Universidade de S\~ao Paulo,\\  S\~ao Carlos, 13560-970,  Brazil
	}

 \newtheorem{thm}{Theorem}[section]
 \newtheorem{cor}[thm]{Corollary}
 \newtheorem{lem}[thm]{Lemma}
 \newtheorem{prop}[thm]{Proposition}
 \theoremstyle{definition}
 
 \theoremstyle{remark}
 \newtheorem{rem}[thm]{Remark}
 
 \numberwithin{equation}{section}
 
 \usepackage{empheq,amsmath}

\usepackage{graphicx,color}

\begin{document}

\maketitle

\begin{abstract}
It is considered a semilinear elliptic partial differential equation in $\mathbb{R}^N$ with a potential that may vanish at infinity and a nonlinear term with subcritical growth.  A positive solution is proved to exist depending on the interplay between the decay of the potential at infinity and the behavior of the nonlinear term at the origin. The proof is based on a penalization argument, variational methods, and  $L^\infty$ estimates. Those estimates allow dealing with settings where the nonlinear source may have supercritical, critical, or subcritical behavior near the origin. Results that provide the existence of multiple and infinitely many solutions when the nonlinear term is odd are also established.
\end{abstract}

\section{Introduction}\label{sec1}

The solutions of the semilinear elliptic equation
\begin{equation}\label{eq1}
-\Delta u + V(x)u = f(x,u),\quad x\in \mathbb{R}^N, \quad N\geq 3,
\end{equation}
play a pivotal role in the search for certain kinds of solitary waves of Klein-Gordon or Schr\"{o}dinger  nonlinear type equations (see \cite{B-L1,B-L2,S}).

The question of the existence of solutions for (\ref{eq1}) has been intensively studied by many researchers under a variety of conditions on $V$ and $f$.  In \cite{FW}, Floer and Weinstein used a Lyapunov-Schmidt type reduction to 
prove the existence of a positive solution of  (\ref{eq1}) for the case $N=1$, with $f(x,s) = s^3$  and $V$  a globally bounded potential having a nondegenerate critical point and 
satisfying $\inf_{\mathbb{R}^N} V >0$.  In \cite{Oh1,Oh2}, Oh generalized the result in \cite{FW} for $N\geq 3$ and  $f(x,s) = \vert s\vert^{p-2}s$, with $2< p < 2^* := 2N/(N-2)$.

The use of mountain pass arguments  to study (\ref{eq1}) goes back at least as far as     \cite{DN} and  \cite{R}.   In \cite{DN}, Ding and Ni obtained several results concerning the existence of positive solutions for the  equations of the type (\ref{eq1}) with $V(x)\geq 0$. Under hypotheses on $f$ coherent with the prototype nonlinearity $Q(x)s^{p-1}$, $2<p< 2^*$, the main results in \cite{DN} establish the existence of positive solutions  for (\ref{eq1}),  decaying uniformly to zero as  $\vert x\vert\to \infty$, when $Q(x)$ is neither radial nor small at infinity.   

In \cite{R}, Rabinowitz, among other results,  employed  a mountain-pass type argument to find a ground state solution for  (\ref{eq1}), supposing $\inf_{\mathbb{R}^N} V >0$   and $f$  superlinear and subcritical at infinity, as long as $\liminf_{\vert x\vert\to \infty} V(x) > 0$  is sufficiently large (see  Theorem 1.7 and Remark 2.21 in  \cite{R}).

 A further step in the study of such problems was made by del Pino and Felmer in \cite{DF}  that, under the hypothesis $\liminf_{\vert x\vert\to \infty} V(x) > 0$,  considered the problem on a local setting by assuming that the potential satisfies $\inf_\Omega V < \inf_{\partial\Omega}V$ on a nonempty bounded open set $\Omega \subset \mathbb{R}^N$. The technique of del Pino and Felmer relies on the study of a  functional  associated with a version of the original problem obtained by a  penalization of the nonlinear term outside $\Omega$. This penalization  allows us to overcome the lack of compactness and to prove the existence of a mountain pass critical level.  Appropriate $L^\infty$ - estimate for the solution
 of the penalized problem implies that it is actually a solution to the original one.
 
 In the last two decades many  studies have focused on potentials that may vanish  at infinity,  that is,  such that $\liminf_{\vert x\vert\to \infty} V(x) = 0$ (we refer the reader to the articles
 \cite{Alves-Souto, Ambrosetti-Malchiodi-Ruiz,  Ambrosetti-Wang, Ambrosetti-Felli-Malchiodi, Azzollini-Pomponio, BPR, BR, BGM, C-M2,  GM} and references therein). It is worth mentioning that in this case the  natural space  to deal with equation (\ref{eq1}) via critical point theory is $D^{1,2}(\mathbb{R}^N)$ which  is only embedded in $L^{2^*}(\mathbb{R}^N)$.   Consequently conditions as  $\vert f(x,s)\vert\leq  C\vert s\vert^{p-1}$,  $p\neq 2^*$, do not guarantee in general  that the related energy functional is well defined.  
 
 Among several methods that have been used to deal with problems with vanishing potentials, we would like to mention a clever adaptation of the penalization technique that has   been introduced by Alves and Souto \cite{Alves-Souto}.  We note that, in order to prove  that the solution of the modified problem is a solution of the original problem, those authors adapt some ideas found in \cite{Brezis-Kato} to obtain an $L^\infty$ - estimate for 
the solution of the penalized problem in terms of its $L^{2^*}$ norm. 
 
  In \cite{Alves-Souto} it is assumed that the potential $V:\mathbb{R}^N \to \mathbb{R}$ is a continuous function satisfying the assumptions:
\begin{equation}
{V(x) \geq 0\ \mbox{ for every $x\in \mathbb{R}^N$,}} \label{V12}
\end{equation}
\begin{equation}
\mbox{there are } \Lambda>0 \mbox{ and }  {R > 1} \mbox{ such that }  \displaystyle \frac{1}{R^4}\inf_{\vert x\vert\geq R} \vert x\vert^4 V (x) \geq \Lambda.\label{V_2}
\end{equation}
Furthermore it is required that  $f$ is  an autonomous continuous function that it  is positive on $(0,\infty)$  and satisfies  the  hypothesis:
\begin{gather}
\limsup_{s\to 0^+} \frac{sf(s)}{s^{2^*}}< +\infty. \label{f_1}
\end{gather}

Assuming further that $f$ is superlinear and subcritical at infinity, the main result in \cite{Alves-Souto} provides a constant $\Lambda^*>0$ such that equation (\ref{eq1}) has a positive solution whenever  $\Lambda  \geq \Lambda^*$. We observe that  condition (\ref{f_1}) implies that $f$ has critical or supercritical decay at the origin.



    Motivated by the articles \cite{Alves-Souto, DF, R}, our primary goal in this paper is to present a version of the penalization technique that will enable us to obtain results on the existence and multiplicity of solutions for equation (\ref{eq1}), depending on the decay of the potential at infinity, under versions of hypothesis (\ref{f_1}) that allow the nonlinear term to have supercritical, critical, or subcritical behavior near the origin.  Furthermore we contemplate the possibility of having $f$ nonautonomous and assuming negative values. We also establish results on the existence of multiple and infinitely many solutions when $f$ is odd with respect to the second variable.

We emphasize that the argument  we use for the existence of solutions  uncovers  an interplay between the behavior of the nonlinear term at the origin and the decay of the potential  at infinity.  A key  ingredient  for establishing such  relation  it is a result of  the $L^\infty$-estimate for  the solution of the penalized problem that does not depend on the behavior of the nonlinear term close to the origin.   



{Henceforth in this article we set $\Omega := \{x\in \mathbb{R}^N : V(x) <0\}$ and define $V_\infty := \max\{ V(x) : \vert x - x_0 \vert \leq r_0 \}$, where $r_0 =1$ and $x_0=0$ if  $\Omega = \emptyset$,  while $r_0>0$ and $x_0 \in \mathbb{R}^N$ are such that $B_{r_0}(x_0) \subset \Omega$ if $\Omega \neq \emptyset$.}

In our first result, we suppose  the function $f$ and the potential $V$ satisfy
\begin{itemize}
\item[(${f_1}$)]  $f \in C(\mathbb{R}^N\times \mathbb{R},\mathbb{R})$ and 
there is $q>2$ such that
$$
\limsup_{s\to 0}\left\vert \frac{sf(x,s)}{s^{q}}\right\vert< +\infty, \ \ \mbox{uniformly in}\  \mathbb{R}^N,
$$
\item[($f_2$)] there  are $a_1, a_2>0$ and  $p\in (2,2^*)$ such that
$$
\vert f(x,s)\vert \leq a_1\vert s\vert ^{p-1}  + a_2\ \  \mbox{for every}\  (x,s)\in \mathbb{R}^N\times \mathbb{R},
$$
\item[($f_3$)]  there are $\theta >2$ and $S_0 \geq 0$ such that
$$
sf(x,s) \geq \theta F(x,s) >0\ \ \mbox{for every}\ \vert s\vert\geq S_0 ,\  x\in \mathbb{R}^N,
$$
where $F(x,s) :=  \int_0^s f(x,t)\,dt$,

\item[$(V_1)$]  $V\in C(\mathbb{R}^N, \mathbb{R})$ and either $V\geq 0$ in $\mathbb{R}^N$ or  $\Omega$ is a nonempty bounded set and 
\[
\inf_{\Omega}V > -S/\vert\Omega\vert^{2/N},
\]
 where $S>0$ is the best constant for the embedding $D^{1,2}(\mathbb{R}^N)$ into $L^{2^*}(\mathbb{R}^N)$,
 
\item [(${V_2}$)] there is {$R> r_0+|x_0|$} such that
\[
{\Lambda(R):= \inf_{\vert x\vert \geq {R}}\vert x\vert^{(N-2)(q-2)} V (x) >0,}
\]
 with  $q>2$ given by  (${f_1}$).
\end{itemize}

We  note  that  ($f_3$) is the  famous  Ambrosetti-Rabinowitz's superlinear condition which has been  introduced  in the seminal article \cite{Ambrosetti-Rabinowitz} to study via minimax methods  semilinear elliptic problems on bounded domains. Condition ($f_3$) with $S_0=0$ has been assumed in the articles \cite{Alves-Souto, DF, R}. Such  hypothesis with $S_0>0$ allows the nonlinear term $f$  to assume negative values.    

It is important   to observe that hypotheses ($f_1$) and ($V_2$) provide a clear relation between the behavior of $f$ near the origin and the lower bound for the  decay of the potential at infinity.  Finally we note that hypothesis  ($f_2$) implies that $f$ has a subcritical growth at infinity and that, under  condition ($V_1$),  the potential $V$ may  assume negative values.



We may now state our first result on the existence of a positive solution for (\ref{eq1}).

\begin{thm}\label{thm}
Suppose $V$ satisfies $({V_1})$-$({V_2})$  and  $f$  satisfies $({f_1})$-$({f_3})$.  Then there is $\Lambda^* >0$ such that  (\ref{eq1}) possesses a positive solution  {provided  $\Lambda(R) \geq \Lambda^*$.}
\end{thm}

Theorem \ref{thm} may be seen as a complement of the result established in \cite{Alves-Souto} because under its hypotheses   $f$ may assume negative values and it does not have  necessarily a critical or supercritical behavior near the origin. In other words, if $q \in (2, 2^*)$,  we obtain a result that is not covered in \cite{Alves-Souto}. Moreover if $q=2^*$, hypothesis $(f_1)$ is exactly 
(\ref{f_1}) and  $({V_2})$  allows the same decrease at infinity for $V$ as the one considered in \cite{Alves-Souto}. If $q>2^*$, Theorem \ref{thm} improves the result of [1] in relation to the behavior (decay) of V at infinity.

{We emphasize that the constant  $\Lambda ^*$ provided  in Theorem \ref{thm} depends on the radius $R>0$ given in condition $(V_2)$. In particular, when the condition $(f_3)$ holds with $S_0=0$, we may find  $\widetilde{\Lambda}^*> 0$, independent of $R$,  such that (\ref{eq1})  has a positive solution whenever  
$\widetilde{\Lambda}(R):= \Lambda(R)/R^{(N-2)(q-2)} \geq \widetilde{\Lambda}^*> 0$  (see Theorem \ref{resultalso} in Section \ref{proofsthm12}).   Observe that  when $q= 2^*$ the corresponding estimate is 
 precisely the hypothesis (\ref{V_2}) assumed in \cite{Alves-Souto}.}




To reinforce the interplay between the behavior of the nonlinear term at the origin and the decay of the potential, evidenced by  conditions $(f_1)$ and  $({V_2})$, we present a result in which  the function $f$ approaches zero rapidly at the origin: supposing  that $f$ and $V$ satisfy
\begin{itemize}
\item[($\widehat{f_1}$)] there are constants $q, a >0$ such that
$$
\displaystyle\limsup_{s\to 0} \vert f(x,s)\vert  e^{(a/\vert s\vert^{q})}< +\infty\ \ \mbox{uniformly in}\  \mathbb{R}^N,
$$
\end{itemize}
\begin{itemize}
\item [({${V_3}$})] there are constants  $\mu >0$ and {$R> r_0+|x_0|$}
such that
\[
{\Lambda(R,\mu):=\inf_{ \vert x\vert \geq {R}}e^{\mu \vert x\vert^{(N-2)q}}V (x) >0,} 
\]
with $q$ given by  ($\widehat{f_1}$),
\end{itemize}
 we may state

\begin{thm}\label{thm2}
Suppose $V$ satisfies $({V_1})$ and $({{V_3}})$,  and  $f$  satisfies  $(\widehat{f_1})$,  $({f_2})$, and $({f_3})$.  Then there are constants $\mu^*, \Lambda^*>0$  such that  (\ref{eq1}) possesses a positive solution {provided}  $0< \mu \leq \mu^*$ and {$\Lambda(R,\mu) \geq \Lambda^*$}.
\end{thm}



We  remark that under the hypotheses of Theorems \ref{thm} and \ref{thm2} we may actually find solutions $u^+$ and $u^-$ of (\ref{eq1}) with $u^+>0$ and $u^-<0$ in $\mathbb{R}^N$. Furthermore, if we suppose $f$ is odd with respect to the second variable, we may derive the existence of multiple pairs of solutions for (\ref{eq1}). More specifically, assuming

\begin{itemize}
\item[(${f_4}$)] $f(x,-s) = -f(x,s)$  for every $(x,s)\in \mathbb{R}^N\times \mathbb{R}$,
\end{itemize}
we may use our version of the penalization technique and  a minimax critical point theorem for functionals with symmetry due to Bartolo et al. \cite{BBF} (see Theorem \ref{simmpt} in Section 5) to obtain:

\begin{thm}\label{thm3}
Suppose $V$ satisfies $({V_1})$-$({V_2})$   and  $f$  satisfies $({f_1})$-$({f_4})$.
Then, given $l\in \mathbb{N}$,  there is $\Lambda_l^* >0$ such that  (\ref{eq1}) possesses $l$ pairs of nontrivial solutions {provided  $\Lambda(R) \geq \Lambda_l^*$.}
\end{thm}

As a  consequence of Theorem \ref{thm3}, we may consider a setting where equation (\ref{eq1}) may have infinitely many pairs of nontrivial solutions without supposing that the potential
is coercive in $\mathbb{R}^N$. For example, assuming the following version of $({V_2})$:
\begin{itemize}
\item [(${{V_4}}$)] {there is a sequences $(R_j) \subset (0,\infty)$}   such that
\[
{\Lambda_j(R_j) := \inf_{ \vert x\vert \geq {R_j}}\vert x\vert^{(N-2)(q-2)} V (x) >0,}
\]
with $q$ given by  ($f_1$), 
\end{itemize}
we may state:
\begin{prop}\label{prop4}
Suppose $V$ satisfies $({V_1})$ and $({{V_4}})$   and  $f$  satisfies $({f_1})$, $(f_2)$, $({f_4})$ and $(f_3)$ with $S_0 =0$.
Then equation   (\ref{eq1}) possesses infinitely many pairs of nontrivial solutions provided
\begin{equation}\label{condprop4}
\limsup_{j\to \infty}\frac{{\Lambda_j(R_j)}}{R_j^{(N-2)(q-2)}} = \infty.
\end{equation}

\end{prop}

 We  remark that hypotheses  (${{V_4}}$) and (\ref{condprop4}) do not imply that the potential $V$ is coercive. Actually, given any $\alpha\geq 0$ it is possible to find a potential satisfying those conditions and $\liminf_{\vert x\vert\to \infty} V(x) = \alpha$. We also note that, under those assumptions, we may always assume that   $R_j \to \infty$, as $j\to \infty$, and, consequently  that 
 {$R_j > r_0+|x_0|$}  for every $j\in\mathbb{N}$.


It is worthwhile mentioning that Proposition \ref{prop4} holds under hypothesis $(f_3)$ with $S_0>0$ and a version of $(\ref{condprop4})$. Furthermore, it is possible to derive the existence of infinitely many pairs of nontrivial solutions for (\ref{eq1}) when $f$ and $V$ satisfy  ($\widehat{f_1}$) and a version of $({V_3})$, respectively (see statements and remark in Section 5).



The article is organized as follows: In Section 2, we introduce the version of the penalization argument used for proving our results and  we establish the existence of a mountain pass solution for  the penalized problem.  Section 3 is devoted to proving an estimate for the $L^\infty$ norm for the solutions to the modified problem in terms of its $L^{2^*}$ norm. In section 3, we also obtain a result on the decay of the solution of the penalized problem at infinity. In Section 4, we present the proofs of Theorems \ref{thm} and \ref{thm2} and the corresponding results when $(f_3)$ holds with $S_0=0$. The proofs of  Theorem \ref{thm3} and Proposition \ref{prop4} are presented in Section 5.

\section{Preliminaries}

Let $E$ be the subspace of $D^{1,2}(\mathbb{R}^N)$ defined by
\[
E=  \left\{
 u\in D^{1,2}(\mathbb{R}^N)  :  \int_{\mathbb{R}^N}V(x)u^2\,dx < \infty \right\}.
\]

We claim that, under the hypothesis (${V_1}$),
\[
\| u \|:=  \left[\int_{\mathbb{R}^N}\left(\vert\nabla u\vert^2 + V(x)u^2\right)\,dx\right]^{1/2}
\]
is a norm in $E$ and $E$ is continuously embedded in $D^{1,2}(\mathbb{R}^N)$.

Since the claim is trivially verified if $V\geq 0$ in $\mathbb{R}^N$, in order to verify the claim it suffices to suppose that $\Omega\neq\emptyset$: given $u\in D^{1,2}(\mathbb{R}^N)$,  we may use H\"{o}lder's inequality, with $r=2^*/2$ and $r'=N/2$,
and the estimate $\vert u\vert^2_{L^{2^*}(\Omega)}\leq S^{-1}\int_{\mathbb{R}^N}\vert\nabla u\vert^2\, dx$ to obtain
\begin{equation}\label{4}
\int_{\Omega}u^2\,dx \leq \frac{\vert\Omega\vert^{2/N}}{S}\int_{\mathbb{R}^N}\vert\nabla u\vert^2\, dx.
\end{equation}
From (${V_1}$),  there is $\alpha >0$ such that
\begin{equation}\label{4,5}
\inf_{x\in \Omega} V(x) \geq -\alpha > -\frac{S}{\vert\Omega\vert^{2/N}}.
\end{equation}
Thus we may invoke   (\ref{4}) to derive
\begin{equation}\label{E-D12}
 \int_{\mathbb{R}^N}\left(\vert\nabla u\vert^2 + V(x)u^2\right)\,dx \geq \left( 1 - \frac{\alpha \vert\Omega\vert^{2/N}}{S}\right)\int_{\mathbb{R}^N}\vert\nabla u\vert^2\, dx.
\end{equation}
 Observing  that $(1- {\alpha \vert\Omega\vert^{2/N}}S^{-1}) >0$, this concludes the proof of the above claim.

Henceforth in this paper we take $\alpha =0$ and $\Omega = \emptyset$ whenever $V\geq 0$ in $\mathbb{R}^N$. Note that in this setting  the above estimates are satisfied for those values of  $\alpha$ and $\Omega$.

 By solution of (\ref{eq1}), we mean a solution in the sense of distribution, namely a function $u\in D^{1,2}(\mathbb{R}^N)$ such that
\[
\int_{\mathbb{R}^N}\left(\nabla u\nabla\phi + V(x)u\phi\right)\, dx - \int_{\mathbb{R}^N}f(x, u)\phi\, dx =0
\]
for every $\phi\in C^\infty_0(\mathbb{R}^N)$, the space of $C^\infty$-functions with compact support.  In fact, an inspection of the proof of our  results will reveal that the weak solution obtained satisfies the above identity for every $\phi\in E$.

    Since in Theorems \ref{thm} and \ref{thm2} we intend to prove the existence of positive solutions, we let $f(x,s)=0$ for every $(x,s)\in \mathbb{R}^N\times (-\infty, 0]$.

Next, in order to deal with the fact  that in our setting $f$ may assume negative values, we introduce a version of the penalization argument employed in \cite{Alves-Souto}:
for $\theta >2$ and $R>0$ given by conditions ($f_3$) and (${V_2}$) respectively, take
$k = 2\theta/(\theta-2)$ and consider  for every $(x,s)\in \mathbb{R}^N\times (0,\infty)$
\[
\tilde{f}(x,s) = \left\{ \begin{array}{rl}

-\frac{1}{k}V(x)s & \mbox{if $kf(x,s) < - V(x)s$},\\
f(x,s) & \mbox{if $-V(x)s\leq kf(x,s) \leq V(x)s$},\\
\frac{1}{k}V(x)s & \mbox{if $kf(x,s) > V(x)s$}.
\end{array}
\right.
\]
Furthermore set $\tilde{f}(x,s) = 0$, for every  $(x,s)\in \mathbb{R}^N\times (-\infty, 0]$,  and  define
\begin{equation}\label{defg}
g(x,s) = \left\{ \begin{array}{ll}
f(x,s) & \mbox{for  $(x,s)\in \mathbb{R}^N\times\mathbb{R}$, $\vert x\vert \leq R$},\\
\tilde{f}(x,s) &  \mbox{for $(x,s)\in \mathbb{R}^N\times\mathbb{R}$, $\vert x\vert > R$}.
\end{array}
\right.
\end{equation}
Observe that $g$ is a Carath\'eodory function satisfying
\begin{equation}
\begin{cases}
g(x,s) = 0\ \mbox{for}\   (x,s)\in \mathbb{R}^N\times (-\infty, 0],\\
g(x,s) = f(x,s) \ \mbox{for}\   (x,s)\in \mathbb{R}^N\times \mathbb{R},\ \vert x\vert \leq R,\label{1} \\
\vert g(x,s)\vert \leq \vert f(x,s)\vert \ \mbox{for}\  (x,s)\in \mathbb{R}^N\times \mathbb{R}, \\
\vert g(x,s)\vert \leq \frac{1}{k}V(x)\vert s\vert \ \mbox{for}\   (x,s)\in \mathbb{R}^N\times \mathbb{R}, \ \vert x\vert > R,\\
\end{cases}
\end{equation}
and
\begin{equation}
\begin{cases}
G(x,s) = F(x,s)\ \mbox{for}\   (x,s)\in \mathbb{R}^N\times \mathbb{R},  \  \vert x\vert \leq R,\\
G(x,s) \leq \frac{1}{2k}V(x)s^2\ \mbox{for}\   (x,s)\in \mathbb{R}^N\times \mathbb{R}, \   \vert x\vert > R,\label{2}
\end{cases}
\end{equation}
where $G(x,s) := \int_0^sg(x,t)\, dt$. The auxiliary problem that we associate with (\ref{eq1}) is the following:
\begin{equation}\label{AP}
   \begin{cases}
 -\Delta u + V(x)u = g(x,u), \quad  x\in \mathbb{R}^N,\\ 
   \quad u\in E.
    \end{cases}
   \end{equation}
We observe that any positive solution $u$ of (\ref{AP}) that satisfies $\vert f(x,u)\vert\leq V(x)u/k$ for $\vert x\vert \geq R$ is actually a solution of (\ref{eq1}). 

 As a consequence of  (\ref{1}) and (\ref{2}), the functional
\[
J(u) = \frac12\int_{\mathbb{R}^N}\left(\vert\nabla u\vert^2 + V(x)u^2\right)\, dx - \int_{\mathbb{R}^N}G(x,u)\, dx
\]
is well defined and of class $C^1$ in $(E, \| \cdot\|)$. Moreover
\begin{equation}\label{ws}
J'(u)v =  \int_{\mathbb{R}^N} (\nabla u\nabla v + V(x)u v)\, dx - \int_{\mathbb{R}^N}g(x,u)v dx\   \mbox{for every  $u, v\in E$}.
\end{equation}
Thus, any critical point of $J$ is a weak solution of (\ref{AP}).

\begin{prop}[mountain pass geometry]\label{mpg}  
Suppose $V$ satisfies $(V_1)$-$(V_2)$  and $f$ satisfies $({f_1})$-$({f_3})$. Then
\begin{enumerate}
\item there exist constants $\beta, \rho >0$ such that $J(u) \geq \beta$ for every $u\in E$ such that $\| u\|= \rho$;
\item there exists $e\in E$, $\| e\| > \rho$, such that $J(e) <0$.
\end{enumerate}
\end{prop}

\begin{proof} 
The proof of item $1$ is standard and follows well-known arguments. For the reader's convenience, we give a proof for the case $\Omega \neq \emptyset$: By $({V_1})$ and  $({V_2})$,  $\Omega \subset B_R(0)$ and $V(x)>0$ for every $\vert x\vert \geq R$. Consequently, from  (\ref{4,5}) and  (\ref{2}), we have
\begin{align}\label{3}
\lefteqn{J(u) = \frac{1}{2}\int_{\mathbb{R}^N}\left(\vert\nabla u\vert^2 + V(x)u^2\right)\, dx - \int_{B_R(0)}\! F(x,u)dx - \int_{\mathbb{R}^N \setminus B_R(0)}\! G(x,u)dx} \nonumber \\
&   \geq \frac{1}{2}\int_{\mathbb{R}^N}\left(\vert\nabla u\vert^2 + V(x)u^2\right)\,dx   -  \int_{B_R(0)}\! F(x,u)dx  -\frac{1}{2k} \int_{\mathbb{R}^N \setminus \Omega}\! V(x)u^2 dx \nonumber \\
&  {\geq} \frac{1}{2}\left[1- \frac{\alpha \vert\Omega\vert^{2/N}}{S}\right]\int_{\mathbb{R}^N}\vert\! \nabla u\vert^2\,dx +  \frac{(k-1)}{2k}\int_{\mathbb{R}^N \setminus \Omega}\! V(x)u^2 dx -    \int_{B_R(0)}\! F(x,u) dx \nonumber \\
& \geq   d_1 \int_{\mathbb{R}^N}(\vert\nabla u\vert^2\,dx +  V(x)u^2)\,dx  -    \int_{B_R(0)} F(x,u) dx,
\end{align}
where $d_1:= \min\{(1- {\alpha \vert\Omega\vert^{2/N}S^{-1}})/2, (k-1)/(2k)\}$.

Given $\epsilon >0$, combining $({f_1})$ with $(f_2)$,  we find a constant $C(\epsilon)>0$ such that
\[
\vert F(x,s)\vert\leq \epsilon\vert s\vert^{2} + C(\epsilon)\vert s\vert^{2^*}\ \mbox{for every} \ (x,s)\in \mathbb{R}^N\times \mathbb{R}.
\]
Thus, there are positive constants $d_2=d_2(R)$ and $d_3=d_3(\epsilon)$ such that
\begin{equation}\label{aux2}
 \int_{B_R(0)} F(x,u)\, dx  \leq \epsilon d_2\| u\|^2 + d_3\| u\|^{2^*} \ \mbox{for every}\  u\in E.
 \end{equation}
 The above estimates and  (\ref{3}) give us
\[
J(u)  \geq  d_1\|u\|^2\ -  \epsilon d_2\|u\|^2 + d_3\|u\|^{2^*} \ \mbox{for every}\  u\in E.
\]
Using the above estimate and taking $\epsilon >0$ sufficiently small,
we  conclude the verification of item $1$ by  finding appropriated values of  $\beta, \rho >0$, 

For proving the item $2$, {recalling the definition of $V_\infty$ from Section \ref{sec1},}
 taking a nonnegative function  $\phi \in E\setminus \{0\}$ such that $\rm{supp}(\phi) \subset B_{r_0}(x_0)$, we obtain
\begin{equation}\label{6}
J(t\phi) \leq \frac{t^2}{2}\int_{B_{r_0}(x_0)}(\vert\nabla \phi\vert^2 + V_\infty \phi^2)\, dx - \int_{B_{r_0}(x_0)}F(x,t\phi)\, dx \quad \forall t\geq 0.
\end{equation}
By $(f_2)$ and  $(f_3)$, there are constants $C_1, C_2>0$, depending on $r_0$,  such that
\begin{equation}\label{7}
F(x,s) \geq C_1s^\theta - C_2\ \mbox{for every}\  (x,s)\in B_{r_0}(x_0)\times [0,\infty).
\end{equation}
From (\ref{6})-(\ref{7}), we have
\begin{equation}\label{estmpl}
J(t\phi) \leq \frac{t^2}{2}\int_{B_{r_0}(x_0)}(\vert\nabla \phi\vert^2 + V_\infty\phi^2)\, dx - C_1t^{\theta}\int_{B_{r_0}(x_0)}\vert\phi\vert^{\theta}\, dx + C_2 \vert B_{r_0}(x_0)\vert
\end{equation}
for every $t\geq 0$. Since $\theta >2$, we have $J(t\phi) \to -\infty$ as $t\to +\infty$. Hence, taking $e= t\phi$, with $t>0$ sufficiently large, we have that  $\| e\| > \rho$ and $J(e) <0$. This concludes the verification of item 2. The proof of Proposition \ref{mpg} is complete. \end{proof}

Recalling that $J$ satisfies the Palais-Smale condition \cite{Ambrosetti-Rabinowitz, R} if every sequence $(u_n) \subset E$ such that  $(J(u_n)) \subset \mathbb{R}$ is bounded and $J'(u_n)\to 0$, as $n \to \infty$, has a strongly convergent subsequence, we may state:

\begin{lem}\label{ps}
Suppose $V$ satisfies $(V_1)$-$(V_2)$  and $f$ satisfies $({f_1})$-$({f_3})$.
Then $J$ satisfies the Palais-Smale condition.
 \end{lem}

\begin{proof}
Let $(u_n) \subset E$ be a sequence such that  $(J(u_n)) \subset \mathbb{R}$ is bounded and $J'(u_n)\to 0$, as $n \to \infty$. First of all we shall verify that   $(u_n)$ is bounded in $(E,\|\cdot\|)$.
Using (\ref{1}), (\ref{2}), we have
\begin{align}\label{9}
J(u_n) & - \frac{1}{\theta}J'(u_n)u_n
= \frac{(\theta - 2)}{2\theta}\| u_n\|^2
+ \int_{\mathbb{R}^N} \left(\frac{1}{\theta}g(x,u_n)u_n - G(x,u_n)\right)\, dx \nonumber\\
&\geq \frac{(\theta - 2)}{2\theta}\|u_n\|^2  + \int_{B_R(0)}\left(\frac{1}{\theta}f(x,u_n)u_n - F(x,u_n)\right)\, dx \nonumber\\
& \quad -  \frac{(\theta +2)}{2\theta k}\int_{\mathbb{R}^N\setminus B_R(0)}V(x)u_n^2\, dx
\end{align}
Next, we invoke $({f_1})$ and $({f_3})$ to find a positive constant {$C=C(R,S_0)$} such that
\begin{equation}\label{estamrab}
\frac{1}{\theta}f(x,s)s - F(x,s) \geq -C\ \mbox{for every}\ (x,s)\in \mathbb{R}^N\times \mathbb{R}.
\end{equation}
Hence (\ref{9}) and our choice of $k$ provide
\begin{align}
J(u_n)  -  \frac{1}{\theta}J'(u_n)u_n &  \geq  \frac{(\theta - 2)}{2\theta}\int_{\mathbb{R}^N}\vert\nabla u_n\vert^2\,dx  + \frac{(\theta - 2)^2}{4\theta^2}\int_{\mathbb{R}^N\setminus B_R(0))} V(x)u_n^2\, dx \nonumber \\
&\quad  + \frac{(\theta - 2)}{2\theta}\int_{B_R(0)} V(x)u_n^2\, dx  - C\vert B_R(0)\vert. \nonumber
\end{align}
Consequently, using  $\theta>2$ and the estimates (\ref{4}) and (\ref{4,5}), with $\alpha =0$ if $\Omega= \emptyset$, we get
\[
J(u_n)  -  \frac{1}{\theta}J'(u_n)u_n   \geq  K\|u_n\|^2  - C\vert B_R(0)\vert,
\]
where
\begin{equation}\label{aux_K}
K: = \frac{(\theta - 2)^2}{4\theta^2}\left(1 - \frac{\alpha \vert\Omega\vert^{2/N}}{S}\right).
\end{equation}
Using that $(J(u_n))\subset \mathbb{R}$ is bounded  and that $J'(u_n)\to 0$, as $n \to \infty$, we conclude that the sequence $(u_n)$ is  bounded in
$(E,\| \cdot\|)$.    Consequently there is $L>0$ such that $\|u_n\| \leq L$ for every $n$. Furthermore, using that $E$ is continuously embedded in $D^{1,2}(\mathbb{R}^N)$ and the Sobolev embedding theorem, by taking a subsequence if necessary, we may suppose there is $u\in E$ such that
\begin{equation}\label{conv}
\begin{cases}
u_n \rightharpoonup u   \mbox{ weakly in $E$},\\
u_n \to u    \mbox{ strongly in $L^{\sigma}$, $\sigma \in [1, 2^*)$, on bounded subsets of $\mathbb{R}^N$},\\
u_n(x) \to u(x)   \mbox{ for almost every $x\in \mathbb{R}^N$}.
\end{cases}
\end{equation}
A standard argument shows that $u$ is a critical point of $J$. Furthermore $u^- = \min(u,0) =0$. Indeed, by (\ref{1}), (\ref{9}) and (\ref{conv}), we have
\[
\|u_n^-\|^2 =  \int_{\mathbb{R}^N} \left(\vert\nabla u_n^-\vert^2 + V(x)(u_n^-)^2\right)\, dx = J'(u_n)u_n^- = o_n(1).
\]
Hence, $u_n^- \to 0$ strongly in $E$ and, consequently,  $u^- = 0$. 
Setting $v_n := u_n - u$, by (\ref{9}) and (\ref{conv}), we have
\begin{align}\label{conv2}
\| v_n\|^2 & = \int_{\mathbb{R}^N}g(x,u_n)v_n\,dx + J'(u_n)v_n -    \int_{\mathbb{R}^N}\left( \nabla u\nabla v_n + V(x)uv_n \right)\,dx \nonumber\\
& = \int_{\mathbb{R}^N}g(x,u_n)v_n\,dx + o_n(1).
\end{align}
Given $\epsilon >0$, we find $r>R$  such that
\begin{equation}\label{conv3}
\int_{\mathbb{R}^N \setminus B_r(0)} V(x)u^2\,dx < \frac{\epsilon^2}{(4k+2)^2L^2}.
\end{equation}
Next, invoking  (\ref{1}), we obtain
\[
\left\vert \int_{\mathbb{R}^N \setminus B_r(0)}\!  g(x,u_n)v_n\, dx\right\vert \leq \frac{1}{k}\int_{\mathbb{R}^N \setminus B_r(0)}\! V(x)u_n^2\,dx + \frac{1}{k}\int_{\mathbb{R}^N \setminus B_r(0)}\! V(x)\vert u_n\vert \vert u\vert\,dx.
\]
The above estimate and (\ref{conv2}) provide
\begin{align}\label{conv4}
\|v_n\|_{D^{1,2}}^2 + \int_{B_r(0)} V(x)v_n^2\,dx + \frac{(k-1)}{k}\int_{\mathbb{R}^N \setminus B_r(0)} V(x)(u_n^2 + u^2)\,dx\nonumber\\
 \leq \int_{B_r(0)}g(x,u_n)v_n\, dx + \frac{(2k+1)}{k}\int_{\mathbb{R}^N \setminus B_r(0)}V(x)\vert u_n\vert\vert u\vert\,dx + o_n(1).
\end{align}
Using that $v_n^2 \leq 2(u_n^2 + u^2)$  in the third integral on the left-hand side of  inequality (\ref{conv4}), we find
\begin{align*}
(k &-1)\|v_n\|^2  \leq 2k\int_{B_r(0)}g(x,u_n)v_n\,dx - (k+1)\int_{\Omega}V(x)v_n^2\,dx \\
& \quad  + (4k+2)\int_{\mathbb{R}^N\setminus B_r(0)}V(x)\vert u_n\vert\vert u\vert\,dx + o_n(1)\\
& \leq 2k\int_{B_r(0)}g(x,u_n)v_n\,dx - (k+1)\int_{\Omega}V(x)v_n^2\,dx \\
& \quad  +  (4k+2)\left[\int_{\mathbb{R}^N\setminus B_r(0)}V(x)u_n^2\,dx\right]^{1/2}\left[\int_{\mathbb{R}^N\setminus B_r(0)}V(x)u^2\,dx\right]^{1/2} + o_n(1).
\end{align*}
Since, by  (\ref{4}) and (\ref{4,5}),
\[
\int_{\mathbb{R}^N\setminus B_r(0)}V(x)u_n^2\,dx \leq \|u_n\|^2 \leq L^2,
\]
we may invoke (\ref{conv}) and (\ref{conv3}) to conclude that $\limsup_{n\to \infty}\|v_n\|^2 < \epsilon$ for every $\epsilon >0$. The fact that $\epsilon >0$ can be chosen arbitrarily small implies that $u_n \to u$ strongly in $E$.
\end{proof}

\begin{rem}\label{rem} We observe that  the decay of $V$ at infinity is not used in the proofs of  Proposition \ref{mpg} and Lemma \ref{ps}. Actually, in those proofs, we have only
 used hypothesis  $({V_1})$  and the fact that $V$  is positive on  $\mathbb{R}^N\setminus B_R(0)$.
\end{rem}

By the Mountain Pass Theorem \cite{Ambrosetti-Rabinowitz}, there is  $u \in E$ such that
\begin{equation}\label{8}
J(u) = c>0\quad \mbox{and}\quad J'(u) = 0,
\end{equation}
where
\begin{equation}\label{level}
c= \inf_{\gamma \in \Gamma}\max_{t\in [0,1]} J(\gamma(t)),
\end{equation}
with $$\Gamma = \{ \gamma \in C([0,1],E) : \gamma(0)=0, \gamma(1) = e\}$$
 for $e$ given by   Proposition \ref{mpg}.

Since $c>0$ and $g(x,s)=0$ for every  $(x,s)\in \mathbb{R}^N\times (-\infty, 0]$, the function $u$ is a nontrivial and nonnegative weak solution of  (\ref{AP}).  Consequently, by the regularity theory and maximum principle, $u$ is positive in $\mathbb{R}^N$. It remains to verify   that $u$ is a solution of (\ref{eq1}).

 We conclude this section with an estimate  for the norm of the solution given by (\ref{8}) that will be of use later to estimate the decay at infinity of the positive weak solutions of (\ref{AP}).  Considering $B_0 := B_{r_0}(x_0)$ and $\phi$ and the constants $C_1$ and $C_2$ given in the proof of Proposition \ref{mpg}, we define
 \begin{equation}\label{supmpl}
d := \sup_{t\geq 0}\left[ \frac{t^2}{2}\int_{B_0}(\vert\nabla \phi\vert^2 + V_\infty\phi^2)dx - C_1t^{\theta}\int_{B_0}\vert\phi\vert^{\theta}dx + C_2 \vert B_0\vert\right]
 \end{equation}

\begin{cor}\label{2kd}
Suppose $u$ is a solution of (\ref{AP}) such that $J'(u)=0$ and $J(u) =c$, with $c>0$ given by (\ref{level}), then
\[
\|u\|^2 \leq K^{-1}\left[d + C\vert B_R(0)\vert\right],
\]
where $C$, $K$  and $d$ are given by (\ref{estamrab}), (\ref{aux_K}) and  (\ref{supmpl}), respectively.
\end{cor}

\begin{proof} First of all we observe that as a direct consequence  (\ref{level}) and (\ref{supmpl}), $c\leq d$.
Furthermore, arguing as in the proof of Lemma \ref{ps},
\[
c= J(u) = J(u)-\frac{1}{\theta}J'(u)u \geq K\|u\|^2 - C\vert B_R(0)\vert.
\]
Consequently
$
\|u\|^2 \leq K^{-1}\left[ c  +C\vert B_R(0)\vert\right] \leq K^{-1}\left[ d + C\vert B_R(0)\vert\right].
$
\end{proof}

\begin{rem} \label{hypalso}
If we suppose $(f_3)$ with $S_0= 0$, the estimate provided by Corollary \ref{2kd} does not depend on  the value of $R$. Indeed, since in this case the constant $C$ given  by  (\ref{estamrab}) is zero, we have  $\|u\|^2 \leq K^{-1}d$ .
\end{rem}

\section{Estimates}

This section is devoted to establishing an estimate for the $L^\infty$ norm of the solutions $u$  in terms of its  $L^{2^*}$ norm and the decay of the solutions. Here we shall consider the problem (\ref{AP})
with $V$ satisfying $({V_1})$ and $g: \mathbb{R}^N\times \mathbb{R} \to \mathbb{R}$ a {Carath\'{e}odory} function satisfying
\begin{itemize}
\item[$({g_1})$]  there are $R> 0$ and $k>1$ such that
$$
\vert g(x,s)\vert \leq \frac{1}{k} V(x) \vert s\vert \ \mbox{ for every}\ s\in \mathbb{R}  \ \mbox{and for a.e.}\  x\in \mathbb{R}^N\setminus B_R (0);
$$
\item[$({g_2})$] there are  $a_1>0$,  $a_2 \geq 0$ and $p\in (2,2^*)$ such that
$$
\vert g(x,s)\vert\leq a_1\vert s\vert^{p-1} + a_2 \  \mbox{ for every}\ s\in \mathbb{R}  \ \mbox{and for a.e.}\  x\in \mathbb{R}^N.
$$
\end{itemize}
Note that as a direct consequence of $({V_1)}$ and $({g_1})$ we have that $\Omega \subset B_R(0)$ whenever $\Omega \neq \emptyset$. 

The next result is a major step in our proofs of the existence of solutions for (\ref{eq1}) because it provides an estimate for the $L^\infty$ norm for the solutions of the auxiliary problem regardless of the behavior of $g$ near the origin.  In our proof of that estimate, we adapt to our setting some of the arguments  used by Br\'{e}zis and Kato \cite{Brezis-Kato} and Alves and Souto \cite{Alves-Souto}.

\begin{lem}\label{infinity}
Suppose $({V_1})$ and $({g_1})$-$({g_2})$ are satisfied. If $u\in E$ is a  solution of problem (\ref{AP}), then $u\in L^\infty(\mathbb{R}^N)$
\begin{equation}\label{estimate}
\vert u\vert_{L^{\infty}(\mathbb{R}^N)} \leq
\left(C_1 \vert u\vert_{L^{2^*}(\mathbb{R}^N)}^{p-2} + C_2 \vert u\vert_{L^{2^*}(\mathbb{R}^N)}^{p-1} + C_3\right)^{\frac{1}{2^* -p}}
\left(1 + \vert u\vert_{L^{2^*}(\mathbb{R}^N)}\right),
\end{equation}
with the constants $C_1, C_2, C_3 \geq 0$ depending on the values of  $a_1$, $a_2$, $\alpha$, $\vert\Omega\vert$, and $p$. Moreover $ C_2 =0$ whenever $a_2 =0$.
\end{lem}

\begin{proof} {For proving Lemma \ref{infinity}, it suffices to verify  that $u^+\in L^\infty(\mathbb{R}^N)$ and the estimate (\ref{estimate}) holds for 
$\vert u^+\vert_{L^{\infty}(\mathbb{R}^N)}$ because to estimate $\vert u^-\vert_{L^{\infty}(\mathbb{R}^N)}$ it is sufficient to observe that  $-u$ is a solution to (\ref{AP}) with $-g(x,-u)$ replacing $g(x,u)$.} 
Moreover, considering that the proof is simpler
when $\Omega = \emptyset$, we shall present a proof of the lemma under the hypothesis $\Omega\neq \emptyset$. Taking $\tau = 2^*/(p-2)$ and {$\sigma = 2^*/2 \tau^\prime$}, where $\tau^\prime$  is the {exponent  conjugate} to $\tau$,
for every  $m\in \mathbb{N}$, $m >a_2$ , we set $v=(u-a_2)^+$ and we define
\begin{align*}
A_m &:= \{ x\in \mathbb{R}^N : \vert v\vert^{\sigma -1} \leq m\}, \\
v_m & := \left\{ \begin{array}{ll}
v\vert v\vert^{2(\sigma -1)} & \mbox{in } A_m,\\
m^2v   & \mbox{in } B_m :=\mathbb{R}^N\setminus A_m,
\end{array}
\right.\\
w_m & := \left\{ \begin{array}{ll}
v\vert v\vert^{(\sigma -1)} & \mbox{in } A_m,\\
mv   & \mbox{in } B_m.
\end{array}
\right.
\end{align*}
A standard verification implies that $v_m \in E$.  Furthermore we  have that $vv_m = w_m^2$ and
\[
\nabla v_m = \left\{ \begin{array}{ll}
(2\sigma-1)\vert v\vert^{2(\sigma -1)}\nabla v & \mbox{in } A_m\\
m^2\nabla v & \mbox{in } B_m
\end{array}
\right.
\]
and
\[
\nabla w_m = \left\{ \begin{array}{ll}
\sigma\vert v\vert^{(\sigma -1)}\nabla v & \mbox{in } A_m\\
m\nabla v & \mbox{in } B_m.
\end{array}
\right.
\]
Hence, considering that $\sigma^2 >2\sigma -1$,
\begin{align}\label{27}
\int_{\mathbb{R}^N} \vert\nabla w_m\vert^2\,dx &= \sigma^2 \int_{A_m}\vert v\vert^{2(\sigma -1)}\vert\nabla v\vert^2\,dx +  m^2 \int_{B_m}\vert\nabla v\vert^2\,dx \nonumber\\
& =\frac{\sigma^2}{2\sigma -1}\int_{\mathbb{R}^N}\nabla v\nabla v_m\,dx
+ m^2\left(1 - \frac{\sigma^2}{2\sigma -1}\right)
\int_{B_m}\vert\nabla v\vert^2\,dx \nonumber\\
& \leq \frac{\sigma^2}{2\sigma -1}\int_{\mathbb{R}^N}\nabla v\nabla v_m\,dx.
\end{align}
Since  $v_m \in E$, $v_m\geq 0$,  and $v_m =0$ on the set $[ u\leq a_2] :=\{x\in \mathbb{R}^N ;\ u(x)\leq a_2\}$,  we may use that  $u$ is a solution of (\ref{AP}) and $({g_1})$ to obtain
\begin{align*}
\int_{\mathbb{R}^N}(\nabla u\nabla v_m  & + V(x)uv_m)\,dx  =   \int_{\mathbb{R}^N}g(x,u)v_m\,dx\\
&\leq  \int_{B_R(0)}\vert g(x,u)\vert v_m\,dx + \frac{1}{k}\int_{\mathbb{R}^N\setminus B_R(0)}V(x)uv_m\,dx.
\end{align*}
Using the definition of $v$ and  that $\Omega\subset  B_R(0)$, we get
\begin{align}
\int_{\mathbb{R}^N}\nabla v\nabla v_m\, dx  + & \int_{\Omega}V(x)uv_m\,dx  \nonumber \\
+ & \frac{(k-1)}{k} \int_{\mathbb{R}^N\setminus\Omega}V(x)uv_m\,dx  \leq
\int_{B_R(0)}\vert g(x,u)\vert v_m\,dx.
\end{align}
Therefore, from $({V_1})$, $({g_2})$, (\ref{4,5}) and (\ref{27})
\begin{equation}\label{28,6}
\int_{\mathbb{R}^N} \vert\nabla w_m\vert^2 \,dx  \leq  \frac{\sigma^2}{2\sigma -1}\left[
\int_{B_R(0)}(a_1\vert u\vert^{p-1} + a_2)v_m \, dx + \alpha\int_{\Omega}u v_m \,dx\right]
\end{equation}
We claim that 
\begin{equation}\label{3.6}
\int_{\mathbb{R}^N} \vert\nabla w_m\vert^2\, dx \leq  \frac{\sigma^2}{2\sigma -1}\left( a_3 \vert v\vert^{2\sigma}_{L^{2\sigma\tau^{\prime}}(\mathbb{R}^N)}+ a_4  \vert v\vert_{L^{2\sigma \tau^{\prime}}(\mathbb{R}^N)}^{2\sigma -1} \right),
\end{equation}
where
\begin{equation}\label{3.71}
a_3 = 2a_1 \vert u\vert_{L^{2^*}(\mathbb{R}^N)}^{p-2} + \alpha \vert\Omega\vert^{\frac{(p-2)}{2^*}};
\end{equation}
and
\begin{equation}\label{3.72}
a_4 = \left\{\begin{array}{l}
2^{p-1} a_2^{-p}(a_1 a_2^{p-2} +1)(1 + \vert u\vert_{L^{2^*}(\mathbb{R}^N)})^{p-1} + \alpha a_2(1 + {\vert}\Omega\vert)\  \mbox{if}\, a_2> 0,\\
0 \  \mbox{if}\, a_2= 0.
\end{array}
\right.
\end{equation}
Indeed, invoking one more time that $v_m= 0$ whenever $u\leq a_2$ and using that $vv_m = w_m^2$, we have
\begin{align*}
\int_{B_R(0)}&(a_1\vert u\vert^{p-1} + a_2)v_m \, dx \\
& \leq 2^{p-1} a_1\int_{\mathbb{R}^N} \vert u\vert^{p-2} w_m^2\, dx  + 2^{p-1} a_2(a_1 a_2^{p-2} + 1)\int_{[u\geq a_2]}v_m \, dx
\end{align*}
Supposing $a_2>0$ and applying H\"{o}lder's inequality, we obtain
\begin{align}\label{3.4}
\int_{B_R(0)}&(a_1\vert u\vert^{p-1} + a_2)v_m \, dx \leq 2^{p-1} a_1\vert u\vert_{L^{2^*}(\mathbb{R}^N)}^{p-2}\left(\int_{\mathbb{R}^N} w_m^{2\tau^\prime}\, dx\right)^{\frac{1}{\tau^\prime}} \\
 & +  2^{p-1}a_2(a_1a_2^{p-2} +1) \vert[u\geq a_2]\vert^{\frac{1}{\tau}}
\left(\int_{[u\geq a_2]} \vert v_m\vert^{\tau^{\prime}}\, dx\right)^{\frac{1}{\tau^\prime}}.\nonumber
\end{align}
Since $\vert w_m\vert\leq \vert v\vert^\sigma$ and $\vert v_m\vert\leq \vert v\vert^{2\sigma -1}$ in $\mathbb{R}^N$,  we have
\begin{align}
\int_{B_R(0)}&
(a_1\vert u\vert^{p-1} + a_2)v_m \, dx \leq 2^{p-1} a_1 \vert u\vert_{L^{2^*}(\mathbb{R}^N)}^{p-2}\vert v\vert_{L^{2\sigma\tau^\prime}(\mathbb{R}^N)}^{2\sigma}\nonumber \\
 & + 2^{p-1}a_2(a_1 a_2^{{p}-2} +1) \vert[u\geq a_2]\vert^{\frac{1}{\tau}}
\left(\int_{[u\geq a_2]}\! \vert v\vert^{(2\sigma -1)\tau^{\prime}}\, dx\right)^{\frac{1}{\tau^\prime}}.\nonumber
\end{align}
Thus, using H\"{o}lder's inequality one more time, we obtain
\begin{align}
\int_{B_R(0)}&(a_1\vert u\vert^{p-1} + a_2)v_m \, dx \leq 2^{p-1}a_1 \vert u\vert_{L^{2^*}(\mathbb{R}^N)}^{p-2}\vert v\vert_{L^{2\sigma\tau^\prime}(\mathbb{R}^N)}^{2\sigma}\nonumber \\
 & + 2^{p-1}a_2(a_1 a_2^{{p}-2} +1) \vert[u\geq a_2]\vert^{(\frac{1}{\tau} + \frac{1}{2\sigma \tau^\prime})}
 \vert v\vert_{L^{2\sigma \tau^{\prime}}(\mathbb{R}^N)}^{2\sigma -1}. \nonumber
\end{align}
Consequently, observing that $a_2^{2^*}\vert[u \geq a_2]\vert\leq \vert  u\vert^{2^*}_{L^{2^*}(\mathbb{R}^N)}$, we may write
\begin{align}
\int_{B_R(0)}&(a_1\vert u\vert^{p-1} + a_2)v_m \, dx \leq 2^{p-1} a_1\vert u\vert_{L^{2^*}(\mathbb{R}^N)}^{p-2}\vert v\vert_{L^{2\sigma\tau^\prime}(\mathbb{R}^N)}^{2\sigma}\nonumber \\
 & +  2^{p-1} a_2^{ -(\frac{2^*}{2\sigma \tau^\prime} + p-1)}(a_1 a_2^{{p}-2} +1)\vert u\vert_{L^{2^*}(\mathbb{R}^N)}^{{(}\frac{2^*}{2\sigma \tau^\prime} + p-2)}
 \vert v\vert_{L^{2\sigma \tau^{\prime}}(\mathbb{R}^N)}^{2\sigma -1}.\nonumber
\end{align}
Since $2^*/2\sigma \tau^\prime= 1$, we obtain
\begin{align}
\int_{B_R(0)}&(a_1\vert u\vert^{p-1} + a_2)v_m \, dx \leq 2^{p-1} a_1\vert u\vert_{L^{2^*}(\mathbb{R}^N)}^{p-2}\vert v\vert_{L^{2\sigma\tau^\prime}(\mathbb{R}^N)}^{2\sigma}\nonumber \\
& +  2^{p-1} a_2^{-p}(a_1 a_2^{{p}-2} +1) (a_2 +1) (1 +\vert u\vert_{L^{2^*}(\mathbb{R}^N)})^{(p-1)}
 \vert v\vert_{L^{2\sigma \tau^{\prime}}(\mathbb{R}^N)}^{2\sigma -1}. \nonumber
\end{align}
Analogously, we obtain
\[
\int_{\Omega} u v_m\, dx \leq \vert\Omega\vert^{\frac{p-2}{2^*}} \vert v\vert^{2\sigma}_{L^{2\sigma\tau^{\prime}}(\mathbb{R}^N)}+ a_2(1 + \vert\Omega\vert)
 \vert v\vert_{L^{2\sigma \tau^{\prime}}(\mathbb{R}^N)}^{2\sigma -1}.
\]
From the above inequalities and (\ref{28,6}), we conclude that the estimate (\ref{3.6}) satisfied if $a_2>0$. A similar argument implies that (\ref{3.6}) {holds} when $ a_2 =0$, and the claim is proved.  From the last claim,  $\vert w_m\vert = \vert v\vert^\sigma$ in $A_m$, the embedding $D^{1,2}(\mathbb{R}^N)\hookrightarrow L^{2^*}(\mathbb{R}^N)$ and $2\sigma -1 >1$,  we obtain
\[
 \left(\int_{A_m} \vert v\vert^{2^*\sigma}\,dx\right)^{2/2^*}
\leq \sigma^2S^{-1}\left( a_3 \vert v\vert^{2\sigma}_{L^{2\sigma\tau^{\prime}}(\mathbb{R}^N)}+ a_4  \vert v\vert_{L^{2\sigma \tau^{\prime}}(\mathbb{R}^N)}^{2\sigma -1} \right).
 \]
Letting $m\to\infty$ and using the monotone convergence theorem, we may write
\[
 \vert v\vert_{L^{2^*\sigma}(\mathbb{R}^N)}^{2\sigma}
\leq \sigma^2 S^{-1}(a_3 + a_4)\left( \vert v\vert^{2\sigma}_{L^{2\sigma\tau^{\prime}}(\mathbb{R}^N)}+  \vert v\vert_{L^{2\sigma \tau^{\prime}}(\mathbb{R}^N)}^{2\sigma -1} \right).
\]
As $\sigma = 2^*/2\tau^\prime$ and replacing $\sigma$ by $\sigma^j$, $j\in \mathbb{N}$, in the above inequality, we obtain
\begin{equation}\label{eq3}
 \vert v\vert_{L^{2^*\sigma^j}(\mathbb{R}^N)}^{2\sigma^j}
\leq \sigma^{2j} S^{-1}(a_3 + a_4)\left( \vert v\vert^{2\sigma^j-1}_{L^{2^*\sigma^{(j-1)}}(\mathbb{R}^N)}+  \vert v\vert_{L^{2^*\sigma^{(j-1)}}(\mathbb{R}^N)}^{2\sigma^j} \right).
\end{equation}
 Using an argument of induction, we may verify that the following inequality holds for every $j\in \mathbb{N}$,
\[
\vert v\vert_{L^{2^*\sigma^j}}\leq \sigma^{\sum_{i=1}^j{i}/{\sigma^i}}\left[2(S^{-1}(a_3 + a_4) +1)\right]^{\frac{1}{2}\sum_{i=1}^j {1}/{\sigma^i}}\left( 1 + \vert v\vert_{L^{2^*}(\mathbb{R}^N)}\right).
\]
Letting  $j\to \infty$ and using that
$
\sum_{i=1}^\infty{i}/{\sigma^{i}} = {\sigma}/{(\sigma -1)^2}$
and $\frac12\sum_{i=1}^\infty1/{\sigma^{i}} = {1}/{2(\sigma -1)}$,
we have
\[
\vert v\vert_{L^{\infty}}\leq \sigma^{\frac{\sigma}{(\sigma-1)^2}}\left[2(S^{-1}(a_3 + a_4) +1)\right]^{\frac{1}{2(\sigma-1)}}\left( 1 + \vert v\vert_{L^{2^*}(\mathbb{R}^N)}\right).
\]
Using that $1/2(\sigma -1) = 1/(2^*-p)$, $a_2\geq 0$, and  $\vert v\vert_{L^{2^*}(\mathbb{R}^N)} \leq \vert u\vert_{L^{2^*}(\mathbb{R}^N)}$, we have
\begin{equation}\label{aux7}
\vert v\vert_{L^{\infty}}\leq 
\left[  \sigma^{\frac{2\sigma}{(\sigma-1)}} 2S^{-1}(a_3 + a_4) + 2\sigma^{\frac{2\sigma}{(\sigma-1)}} + a_2^{2(\sigma-1)}\right]^{\frac{1}{2^* - p}}
(1 + \vert u\vert_{L^{2^*}(\mathbb{R}^N)}).
\end{equation}
{
Since $2(\sigma -1)/(2^*-p) =1$, we have
\begin{equation}\label{aux8}
a_2 < 
\left[  \sigma^{\frac{2\sigma}{(\sigma-1)}} 2S^{-1}(a_3 + a_4) + 2\sigma^{\frac{2\sigma}{(\sigma-1)}} + a_2^{2(\sigma-1)}\right]^{\frac{1}{2^* - p}}
(1 + \vert u\vert_{L^{2^*}(\mathbb{R}^N)}). 
\end{equation}
As $u^+    \leq   (u - a_2)^+ +  a_2  = v + a_2$,  by (\ref{aux7})- (\ref{aux8}), we get
\[
\vert u^+ \vert_{L^{\infty}}\leq 
2\left[  \sigma^{\frac{2\sigma}{(\sigma-1)}} 2S^{-1}(a_3 + a_4) + 2\sigma^{\frac{2\sigma}{(\sigma-1)}} + a_2^{2(\sigma-1)}\right]^{\frac{1}{2^* - p}}
(1 + \vert u\vert_{L^{2^*}(\mathbb{R}^N)}).
\]}
Hence, from (\ref{3.71}) and (\ref{3.72}),
\[
\vert u^+ \vert_{L^{\infty}}\leq
\left[  C_1 \vert u\vert_{L^{2^*}(\mathbb{R}^N)}^{p-2} + C_2 \vert u\vert_{L^{2^*}(\mathbb{R}^N)}^{p-1} + C_3 \right]^{\frac{1}{2^* - p}}
(1 + \vert u\vert_{L^{2^*}(\mathbb{R}^N)}),
\]
where $C_1, C_2$ and $C_3$ are constants depending on the values of $a_1,a_2, p, \alpha, \vert\Omega\vert$ and $k$. Moreover, $C_2 =0$ whenever  $a_2=0$. The proof of Lemma \ref{infinity} is complete.
\end{proof}

\begin{rem}\label{obsCinfty}
We note that the values of the constants $C_1$, $C_2$, and $C_3$ given by Lemma \ref{infinity} do not depend on the value of $R$ in hypothesis $(g_1)$ or on the value of the potential $V$ on $\mathbb{R}^N\setminus \Omega$.
\end{rem}


\begin{lem}\label{decay} Suppose ($V_1$)  and $(g_1)$-$(g_2)$ are satisfied.
If $u\in E$ is  a weak solution of (\ref{AP}), then
\[
\vert u(x)\vert \leq  M\left(\frac{R}{\vert x\vert}\right)^{N-2} \ \mbox{for every} \ x \in \mathbb{R}^N,\   \vert x\vert\geq R,
\]
 where $R>0$ is given by $(g_1)$ and
\[
M= \left(C_1 \vert u\vert_{L^{2^*}(\mathbb{R}^N)}^{p-2} + C_2 \vert u\vert_{L^{2^*}(\mathbb{R}^N)}^{p-1} + C_3\right)^{\frac{1}{2^* -p}}
\left(1 + \vert u\vert_{L^{2^*}(\mathbb{R}^N)}\right),
\]
with  $C_1$ and $C_2$  the constants given  by Lemma \ref{infinity}.
\end{lem}

\begin{proof}

Let $z \in C^{\infty}(\mathbb{R}^N\setminus\{0\})$ be the harmonic function $z(x) = M (R/\vert x\vert)^{N-2}$, for $x\in \mathbb{R}^N$.  Next,  take
\[
w^+(x) = \left\{ \begin{array}{ll}
(u(x) - z(x))^+ & \mbox{if}\  \vert x\vert \geq R,\\
0 & \mbox{if}\ \vert x\vert < R.
\end{array}
\right.
\]
As a direct consequence of Lemma \ref{infinity} $\vert u(x)\vert \leq z(x)$, for every $x$ such that $\vert x\vert=R$. Moreover, since  $\Delta z=0$  in $\mathbb{R}^N\setminus B_R(0)$,  $w^+\in E$, $w^+(x) =0$ for every $\vert x\vert\leq R$, and $w^+\geq 0$,  by $(g_1)$, we have
\begin{align*}
\int_{\mathbb{R}^N}\vert\nabla w^+\vert^2\, dx & =\int_{\mathbb{R}^N}\nabla(u-z)\nabla w^+\, dx\\
& = \int_{\mathbb{R}^N}\nabla u\nabla w^+\, dx - \int_{\mathbb{R}^N}\nabla z\nabla w^+\, dx\\
&=  \int_{\mathbb{R}^N}\left(g(x,u)w^+ - V(x)uw^+\right)\, dx\\
&\leq \int_{\mathbb{R}^N\setminus B_R(0)}V(x)\left( \frac{\vert u\vert}{k} - u\right)w^+\, dx \leq 0,
\end{align*}
where in the last inequality we have used  that $u\geq 0$ whenever $w^+>0$.  
Now, taking
\[
w^-(x) = \left\{ \begin{array}{ll}
(- u(x) - z(x))^+& \mbox{if}\  \vert x\vert \geq R,\\
0 & \mbox{if}\ \vert x\vert < R,
\end{array}
\right.
\]
arguing as above and observing that $u\leq 0$ whenever $w^->0$, we obtain
\begin{align*}
\int_{\mathbb{R}^N}\vert\nabla w^-\vert^2\, dx & =\int_{\mathbb{R}^N}\nabla(-u-z)\nabla w^-\, dx\\
&=  \int_{\mathbb{R}^N}\left(-g(x,u)w^- + V(x)uw^-\right)\, dx\\
&\leq \int_{\mathbb{R}^N\setminus B_R(0)}V(x)\left( \frac{\vert u\vert}{k} + u\right)w^-\, dx \leq 0.
\end{align*}
It follows from the above estimates that  $w ^i\equiv 0$, $i = \pm$, and, consequently, $\vert u\vert\leq z$ on $\mathbb{R}^N\setminus B_R (0)$. The proof of Lemma \ref{decay}
is complete.
\end{proof}

\section{Existence of a positive solution}\label{proofsthm12}
In this section, we present the proofs of Theorems \ref{thm} and  \ref{thm2} and the corresponding results when hypothesis $(f_3)$ holds with $S_0=0$.

\medskip
\noindent{\textbf{Proof of Theorem \ref{thm}}}.
In view  of  Corollary \ref{2kd}, (\ref{E-D12}) and the estimate $\vert u\vert^2_{L^{2^*}(\Omega)}\leq S^{-1}\int_{\mathbb{R}^N}\vert\nabla u\vert^2\, dx$, the equation (\ref{AP}) has a positive solution $u\in E$ satisfying
\begin{equation}\label{4.1}
\vert u\vert_{L^{2^*}(\mathbb{R}^N)} \leq \widehat{C} := \left[ K^{-1}( S - \alpha \vert\Omega\vert^{2/N})^{-1} (d + C \vert B_R(0)\vert)\right]^{1/2},
\end{equation}
with $C$, $K$ and $d$ given by (\ref{estamrab}), (\ref{aux_K}) and (\ref{supmpl}), respectively.
Since $g$ defined by (\ref{defg}) satisfies the hypotheses $(g_1)$ and $(g_2)$, with  $k = 2\theta/(\theta -2)$ and $R$,
$a_1$, $p$ and $\theta$ given by $(V_2)$, $(f_2)$ and $(f_3)$,   in view of (\ref{4.1}) and Lemma \ref{decay}, we have
\begin{equation}\label{4.2}
\vert u(x)\vert \leq  \widehat{M}\left(\frac{R}{\vert x\vert}\right)^{N-2} \ \mbox{for every} \ x \in \mathbb{R}^N,\   \vert x\vert\geq R,
\end{equation}
 where
\begin{equation}\label{4.3}
\widehat{M} :=
\left(C_1\widehat{C}^{p-2} + C_2\widehat{C}^{p-1} + C_3\right)^{1/(2^* -p)}\left( 1 + \widehat{C}\right),
\end{equation}
for constants $C_1$ and $C_2$ given by Lemma \ref{infinity}.
Next, using the hypotheses $(f_1)$ and $(f_2)$, we find  a constant $C>0$ such that
\begin{equation}\label{4.4}
\vert f(x,s)\vert \leq  C\vert s\vert^{q-1}\  \mbox{for every}\ (x,s)\in \mathbb{R}^N\times \mathbb{R}, \ \vert s\vert\leq \widehat{M}.
\end{equation}
Therefore, from (\ref{4.1}) - (\ref{4.4}),
\begin{equation}\label{aux11}
\vert f(x,u(x))\vert\leq  C \widehat{M} ^{(q-2)}\left(\frac{R}{\vert x\vert}\right)^{(N-2)(q-2)}\vert u(x)\vert\ \  \mbox{ for every}\  \vert x\vert\geq R.
\end{equation}
Setting $\Lambda^* =  k C \widehat{M} ^{(q-2)}R^{(N-2)(q-2)}$ and combining  (\ref{aux11}) with $(V_2)$, we obtain
\[
\vert f(x,u(x))\vert \leq \frac{1}{k}V(x)\vert u(x)\vert\ \mbox{ for every}\   \vert x\vert\geq R,
 \]
 whenever $\Lambda(R) \geq \Lambda^*$.  Thus, $\tilde{f}(x,u(x)) = f(x, u(x))$ for every $x\in \mathbb{R}^N$ such that $\vert x\vert \geq R$, and so  $g(x,u(x)) = f(x, u(x))$ for every $x\in \mathbb{R}^N$. Therefore, $u$ is a positive solution of (\ref{eq1}). The proof of Theorem \ref{thm} is complete. \hfill $\Box$

\medskip
As it was observed in the introduction, when $(f_3)$ holds with $S_0=0$, we may provide  a relation between the parameter  in hypothesis $(V_2)$ and the value of $R$ that, in particular, generalizes the result in \cite{Alves-Souto} when $q>2^*$.

\begin{thm} \label{resultalso} 
Suppose $V$ satisfies $({V_1})$ and $({V_2})$,   and  $f$  satisfies  $({f_1})$-$({f_2})$, and $({f_3})$ with $S_0 = 0$.  Then there is $\widetilde{\Lambda}^*>0$ such that  (\ref{eq1}) possesses a positive solution {provided   $\widetilde{\Lambda}(R):= \Lambda(R)/R^{(N-2)(q-2)} \geq \widetilde{\Lambda}^*$.}
\end{thm}

\begin{proof}
  Since $f$ satisfies  $(f_3)$  with $S_0 =0$, we may invoke
Remark \ref{hypalso} to conclude that the constant $\widehat{C}$, given by (\ref{4.1}),
does not depend on the value of $R$. Therefore, from Remark \ref{obsCinfty} and (\ref{4.3}), the estimate
(\ref{4.4}) holds with the constants $C$ and $\widehat{M}$ being independent of the value of $R$. Consequently, supposing that $({V_2})$ holds, the argument used in the proof of Theorem \ref{thm} implies that
(\ref{eq1}) possesses a positive solution whenever $\widetilde{\Lambda}(R):= \Lambda(R)/R^{(N-2)(q-2)} \geq \widetilde{\Lambda}^* := {kC \widehat{M}^{(q-2)}}$.
\end{proof}

\medskip
\noindent{\textbf{Proof of Theorem \ref{thm2}}}.
Since $(\widehat{f_1})$ implies that $(f_1)$ holds, we may invoke Corollary \ref{2kd}, Lemma \ref{decay} and the argument used in the proof of Theorem \ref{thm} to conclude that equation (\ref{AP}) has a positive solution $u\in E$ satisfying the estimate (\ref{4.2}).
Next, we  fix  $0< \widehat{a} < a$. From   $(\widehat{f_1})$ and $(f_2)$ we may find $C>0$ such that
\[
\vert f(x,s)\vert \leq Ce^{-\widehat{a}/ \vert s\vert^q}\vert s\vert\ \mbox{for every}\  (x,s)\in \mathbb{R}^N\times\mathbb{R},\  \vert s\vert\leq \widehat{M},
\]
where  $\widehat{M}$ is given by (\ref{4.3}). Consequently, by (\ref{4.2}),  we may write
\[
\vert f(x,u(x))\vert\leq C e^{-\mu^* \vert x\vert^{(N-2)q}}\vert u(x)\vert \ \mbox{for every}\ x\in \mathbb{R}^N, \ \vert x\vert\geq R,
\]
where $\mu^* := \widehat{a}/ (\widehat{M}R^{(N-2)})^q$. Taking $\Lambda^* = k C$ and using $(V_3)$, we get
\begin{equation}\label{aux13}
\vert f(x,u(x))\vert \leq \frac{1}{k}V(x)\vert u(x)\vert\ \mbox{ for every} \vert x\vert\geq R
 \end{equation}
  whenever  {$0< \mu \leq \mu^*$ and $\Lambda(R,\mu) \geq \Lambda^*$}. This implies  that $u$ is a positive solution of (\ref{eq1}) because 
   (\ref{defg}) and (\ref{aux13})  imply $g(x,u(x)) = f(x,u(x))$ for all
$x\in  \mathbb{R}^N$. The proof of Theorem \ref{thm2} is complete. \hfill $\Box$

\medskip

{Remark similar to Theorem \ref{thm} shows that $\Lambda^*$ in Theorem \ref{thm2} depends on $R$ given in $(V_3)$ if $S_0 >0$ in $(f_3)$.  If $(f_3)$ holds with $S_0=0$ and $V$ satisfies
\begin{itemize}
\item [($V_5$)] {there are constants  $\mu >0$ and {$R> r_0+|x_0|$} such that
\[
\widehat{\Lambda}(\mu,R):= \inf_{\vert x\vert\geq R}e^{\mu \left(\frac{\vert x\vert}{R}\right)^{(N-2)q}}V (x) >0,
\]}
with $q$ given by  ($\widehat{f_1}$),
\end{itemize}
 then we may  find $\widehat{\Lambda}^* > 0$ independent of $R$ such that (\ref{eq1}) has a positive solution provided   $\widehat{\Lambda}(\mu,R) > \widehat{\Lambda}^*$. More precisely,  using the arguments employed in the  proof of Theorem \ref{thm2}, we may also state the  following version of Theorem \ref{resultalso}: }

\begin{thm}\label{thm4.2}
Suppose $V$ satisfies $({V_1})$ and $({V_6})$,  and  $f$  satisfies  $(\widehat{f_1})$,  $({f_2})$, and $({f_3})$ with $S_0=0$.
Then there are constants $\widehat{\Lambda}^*, \widehat{\mu}^* >0$  such that  (\ref{eq1}) possesses a positive solution  {provided 
$0< \mu \leq \widehat{\mu}^*$ and $\widehat{\Lambda}(\mu,R) \geq \widehat{\Lambda}^*$}.
\end{thm}

\section{Existence of multiple solutions}\label{section5}
In this section, we present proofs of Theorem \ref{thm3}  and Proposition \ref{prop4}. Here we also state a version of Proposition \ref{prop4} under hypothesis   $(f_3)$  with 
$S_0>0$ and corresponding results under hypothesis  $(\widehat{f_1})$ and versions of $(V_3)$.

Before proving Theorem \ref{thm3}, we recall an abstract result that provides the existence of $k$ pairs of nontrivial solutions for even functionals defined  in Banach spaces (see 
\cite{BBF}).

\begin{thm}\label{simmpt}
  Suppose $ X= X_1\bigoplus X_2$, with $\dim{X_1} < \infty$,  is a Banach space and that $I\in C^1(X,\mathbb{R})$ is an even functional satisfying $I(0)=0$, the Palais-Smale condition, and
\begin{itemize}
  \item[$(I_1)$] there exist $\beta, \rho >0$ such that $I(u)\geq \beta$  for every $u\in X_2$ such that $\|u\| = \rho$;
  \item[$(I_2)$] there exist a finite dimensional  subspace $W$ of $X$ and $\gamma >0$ such that $I(u) \leq \gamma$ for every $u\in W$.
\end{itemize}
Then, if $\dim{W} >\dim{X_1}$, the functional  $I$ possesses $l$ pairs of nontrivial critical points $\{\pm u_1, \cdots , \pm u_l \}$, where $l = \dim{W} - \dim{X_1}$. Furthermore
$\beta\leq  I(u_i)\leq \gamma$ for every $i\in \{ 1, \cdots , l\}$.
\end{thm}

In order to apply Theorem \ref{simmpt}, we consider $\hat{g} :\mathbb{R}\times \mathbb{R}^N\rightarrow \mathbb{R}$ the odd extension of the function defined by (\ref{defg}) for every $(x,s)\in \mathbb{R}^N \times [0,\infty)$. We note that $\hat{g}$ is a Carath\'{e}odory function satisfying
\begin{equation}
\begin{cases}
\hat{g}(x,s) = f(x,s) \ \mbox{for}\   (x,s)\in \mathbb{R}^N\times \mathbb{R},\ \vert x\vert \leq R, \label{5.1} \\
\vert\hat{g}(x,s)\vert \leq \vert f(x,s)\vert \ \mbox{for}\  (x,s)\in \mathbb{R}^N\times \mathbb{R},\\
\vert\hat{g}(x,s)\vert \leq \frac{1}{k}V(x)\vert s\vert\ \mbox{for}\   (x,s)\in \mathbb{R}^N\times \mathbb{R}, \ \vert x\vert > R.\\
\end{cases}
\end{equation}
Moreover, defining  $\hat{G}(x,s) := \int_0^sg(x,t)\, dt$, we  have
\begin{equation}
\begin{cases}
\hat{G}(x,s) = F(x,s)\ \mbox{for}\   (x,s)\in \mathbb{R}^N\times \mathbb{R},  \  \vert x\vert \leq R, \\
\hat{G}(x,s) \leq \frac{1}{2k}V(x)s^2\ \mbox{for}\   (x,s)\in \mathbb{R}^N\times \mathbb{R}, \   \vert x\vert > R. \label{5.2}
\end{cases}
\end{equation}
Next,  we   consider the symmetric version of the auxiliary equation (\ref{AP}):
\begin{equation}\label{ODDAP}
    \begin{cases}
    -\Delta u + V(x)u = \hat{g}(x,u),\quad  x\in \mathbb{R}^N,\\
     \quad u\in E.
    \end{cases}
    \end{equation}
We also  observe that any solution $u$ of (\ref{ODDAP}) that satisfies $\vert f(x,u)\vert\leq V(x)\vert u\vert/k$, for $\vert x\vert \geq R$, is a solution of (\ref{eq1}).

As in the proofs of Theorems \ref{thm} and \ref{thm2}, the functional   $\hat{J}: E \rightarrow \mathbb{R}$, given by
\[
\hat{J}(u) = \frac12\int_{\mathbb{R}^N}\left(\vert\nabla u\vert^2 + V(x)u^2\right)\, dx - \int_{\mathbb{R}^N}\hat{G}(x,u)\, dx
\]
associated with (\ref{ODDAP}), is  of class $C^1$ in $(E, \| \cdot\|)$ and critical points of $\hat{J}$ are weak solutions of equation (\ref{ODDAP}). Furthermore $\hat{J}$ is an even functional functional and $\hat{J}(0)=0$.

The next result provides a version of Proposition \ref{mpg} for the functional $\hat{J}$.
\begin{prop}\label{oddmpg}
  Suppose $V$ satisfies $(V_1)$-$(V_2)$  and $f$ satisfies $({f_1})$-$({f_3})$. Then
  \begin{enumerate}
\item there exist constants $\beta, \rho >0$ such that $\hat{J}(u) \geq \beta$ for every $u\in E$ such that $\| u\|= \rho$;
\item given $l\in \mathbb{N}$, there exist $\{\phi_1, \cdots , \phi_l\} \subset E$ and $D_l>0$ such that  such that $\hat{J}(u) { \leq } D_l $ for every $u\in W_l := \langle \phi_1, \cdots, \phi_l\rangle$.
\end{enumerate}
\end{prop}
\begin{proof} Item $1$ has been actually proved in Proposition \ref{mpg}. For proving item $2$,  {as in the proof of item (2) of  Proposition \ref{mpg},  
for each  $1\leq i\leq l$  we consider nonnegative functions  $\phi_i \in E\setminus \{0\}$ such that $\rm{supp}(\phi_i) \subset B_0:=B_{r_0}(x_0)$  and  $\rm{supp}(\phi_i)\cap\rm{supp}(\phi_j) = \emptyset$ for  $i\neq j\in \{1,\cdots, l\}$.}  Defining $d_i$, $i\in\{1, \cdots, l\}$,  by
\[
 d_i :=
 \sup_{\tau \in \mathbb{R}}\left[ \frac{\tau^2}{2}\int_{B_0}(\vert\nabla \phi_i\vert^2 + V_\infty\phi_i^2)dx - C_1\vert\tau\vert^{\theta}\int_{B_0}\vert\phi_i\vert^{\theta}dx + C_2 \vert B_0\vert\right],
\]
with $C_1$ and $C_2$ given by (\ref{7}),  we obtain
$J(t\phi_i) \leq d_i <\infty$ for every $t\in \mathbb{R}$. Since $\rm{supp}(\phi_i)\cap\rm{supp}(\phi_j) = \emptyset$ for  $i\neq j\in \{1,\cdots, l\}$, item $2$ is satisfied with $D_l = { \sum_{i=1}^{l}d_i}$. The proposition is proved.
\end{proof}

We are now ready to present the proofs of Theorem \ref{thm3} and Proposition \ref{prop4}:

\medskip
\noindent{\textbf{Proof of Theorem \ref{thm3}}}.
 Observing that  the argument employed in the proof of Lemma \ref{ps} may be easily adapted to verify that $\hat{J}$ satisfies the Palais-Smale condition, we may invoke Theorem  \ref{simmpt} and Proposition \ref{oddmpg} to verify that
$\hat{J}$ possesses $l$ pairs of nontrivial critical points $\{\pm u_1, \cdots , \pm u_l \}$. Moreover
$\beta \leq \hat{J}(u_i) \leq D_l$, and
\begin{equation}\label{estimateui}
\|u_i\| \leq K^{-1} \left( D_l + C \vert B_R(0)\vert\right) \ \mbox{for every} \ i\in \{1,\cdots, l\},
\end{equation}
where $C$, $K$ and $D_l$ are given by (\ref{estamrab}), (\ref{aux_K}) and Proposition \ref{oddmpg}-item 2, respectively. The proof of Theorem \ref{thm3} is concluded by using the estimate provided by Lemma \ref{decay} and arguing as in the proof of Theorem \ref{thm}. \hfill $\Box$

\medskip
\noindent{\textbf{Proof of Proposition \ref{prop4}}}.  Noting that the estimate (\ref{estimateui}) holds with $C=0$ (see Remark \ref{hypalso}), we may argue as in the proof of Theorem
\ref{resultalso} to conclude that under the hypothesis of Theorem \ref{thm3}, with  $f$ satisfying $(f_3)$ with $S_0=0$, we may find $l$ pairs of nontrivial pairs of solutions whenever 
$V$ satisfies  $(V_4)$ with ${\Lambda_j(R_j)} / R_j^{(N-2)(q-2)}$ is sufficiently  large, independently of the value of $R_j$. This fact and the hypothesis  $\limsup_{j\to\infty}{\Lambda_j(R_j)} / R_j^{(N-2)(q-2)}= \infty$ conclude the proof of Proposition \ref{prop4}. \hfill $\Box$

It is not difficult to verify that,  as a direct consequence  Theorem \ref{thm3}, a version of Proposition \ref{prop4} holds under $(f_3)$ without supposing $S_0=0$:

\begin{prop}\label{prop5.3}
Suppose $V$ satisfies $(V_1)$ and $({V_4})$ and $f$ satisfies $(f_1)$-$(f_4)$. Then there exists a sequence $(\Lambda^*_{j,l})\subset (0,\infty)$, $(j,l)\in \mathbb{N}^2$, such that 
(\ref{eq1})  possesses infinitely many pairs of nontrivial solutions whenever $({\Lambda_j(R_j)})$ has a subsequence $({\Lambda_{j_i}(R_{j_i})})$ satisfying ${\Lambda_{j_i}(R_{j_i})}\geq \Lambda^*_{j_i,l_i}$  for every $i\in \mathbb{N}$, with $l_i\to\infty$, as $i\to\infty$.
\end{prop}

\medskip
We also state without proving  a   version of Theorem \ref{thm3} under the hypothesis of Theorem \ref{thm2} and $(f_4)$.

\begin{thm}\label{5.4}
Suppose $V$ satisfies $({V_1})$ and $({V_3})$,  and  $f$  satisfies  $(\widehat{f_1})$, and $({f_2})$-$({f_4})$.
Then, given $l\in \mathbb{N}$ there are constants $\Lambda^*_l, \mu^*_l >0$  such that  (\ref{eq1}) possesses $l$ pairs of nontrivial  solution  {provided 
 $0< \mu \leq \mu^*_l$ and $\Lambda(R,\mu) \geq \Lambda^*_l$}.
\end{thm}

\begin{rem}\label{5.5}
Finally we like to emphasize that, assuming  the condition $(f_4)$,  it is not difficult to  derive versions of Propositions  \ref{prop4} and \ref{prop5.3} under the hypothesis 
$(\widehat{f_1})$ and conditions related to $({V_3})$ and $({V_4})$.
\end{rem}

\subsection*{Acknowledgment}
This work  was financed in part by the CNPq/Brazil. Sergio H. M. Soares wishes to thank the Department of Mathematics of the University of Bras\'{\i}lia, where part of this article was written, for the hospitality.


\begin{thebibliography}{1}

\bibitem{Alves-Souto} C. O. Alves,  M. A. S. Souto,  \emph{Existence of solutions for a class of elliptic equations in $\mathbb{R}^N$ with vanishing potentials.} J. Differential Equations \textbf{252} (2012), 5555--5568.



\bibitem{Ambrosetti-Malchiodi-Ruiz}
A. Ambrosetti, A. Malchiodi, D. Ruiz, \emph{Bound states of nonlinear Schr\"{o}dinger equations with potentials vanishing at infinity.} J. Anal. Math. \textbf{98} (2006), 317--348.

\bibitem{Ambrosetti-Wang} 
A. Ambrosetti,  Z.-Q. Wang,  \emph{Nonlinear Schr\"{o}dinger equations with vanishing and decaying potentials.} Differential Integral Equations \textbf{18} (2005), 1321--1332.

\bibitem{Ambrosetti-Felli-Malchiodi}
A. Ambrosetti, V. Felli, A. Malchiodi, \emph{Ground states of nonlinear Schr\"{o}dinger equations with potentials vanishing at infinity.} J. Eur. Math. Soc. (JEMS) \textbf{7} (2005),  117--144.


\bibitem{Ambrosetti-Rabinowitz}
A. Ambrosetti, P. H.  Rabinowitz,  \emph{Dual variational methods in critical point theory and applications.} J. Funct. Anal. \textbf{14} (1973), 349--381.


\bibitem{Azzollini-Pomponio}
A. Azzollini, A. Pomponio,  \emph{On the Schrödinger equation in $\mathbb{R}^N$ under the effect of a general nonlinear term.} Indiana Univ. Math. J. \textbf{58} (2009), 1361--1378.

\bibitem{BPR} M. Badiale, L. Pisani, S. Rolando, 
\emph{Sum of weighted Lebesgue spaces and nonlinear elliptic equations.}
NoDEA Nonlinear Differential Equations Appl. \textbf{18} (2011), 369--405.

\bibitem{BR} M. Badiale, S. Rolando, 
\emph{Elliptic problems with singular potential and double-power nonlinearity.}
Mediterr. J. Math. \textbf{2} (2005),  417--436.

\bibitem{BBF} 
P. Bartolo, V. Benci, D. Fortunato,  
\emph{Abstract critical point theorems and applications to some nonlinear problems with ``strong'' resonance at infinity.}
Nonlinear Anal. \textbf{7} (1983), 981--1012.



\bibitem{BGM} V. Benci,  C. R. Grisanti, A. M. Micheletti, 
\emph{Existence and non-existence of the ground state solution for the nonlinear Schroedinger equations with $V(\infty)=0$.}
Topol. Methods Nonlinear Anal. \textbf{26} (2005), 203--219.

\bibitem{B-L1} H. Berestycki,   P.-L. Lions, \emph{Nonlinear scalar field equations. I. Existence of a ground state.} Arch. Rational Mech. Anal. \textbf{82} (1983), 313--345.

\bibitem{B-L2} H. Berestycki, P.-L. Lions,  \emph{Nonlinear scalar field equations. II. Existence of infinitely many solutions.} Arch. Rational Mech. Anal. \textbf{82} (1983), 347--375.

\bibitem{Brezis-Kato} H. Br\'{e}zis, T. Kato, \emph{Remarks on the Schr\"{o}dinger operator with regular complex potentials,} J. Math. Pures Appl. \textbf{58} (1979), 137--151.


\bibitem{C-M2} M. Clapp, L. A.  Maia,  \emph{Existence of a positive solution to a nonlinear scalar field equation with zero mass at infinity.} Adv. Nonlinear Stud. \textbf{18} (2018),  745--762.




\bibitem{DF}  M. del Pino, M., P. L. Felmer, 
\emph{Local mountain passes for semilinear elliptic problems in unbounded domains.} Calc. Var. Partial Differential Equations \textbf{4} (1996),  121--137.

\bibitem{DN}  W. Y. Ding, W.-M.  Ni, \emph{On the existence of positive entire solutions of a semilinear elliptic equation.} Arch. Rational Mech. Anal. \textbf{91} (1986),  283--308.



\bibitem{FW} A. Floer, A. Weinstein,  \emph{Nonspreading wave packets for the cubic Schrödinger equation with a bounded potential.}
J. Funct. Anal. \textbf{69} (1986),  397--408.

\bibitem{GM} M. Ghimenti, A. M.  Micheletti,  \emph{Solutions for a nonhomogeneous nonlinear
Schroedinger equation with double-power nonlinearity.} Differ. Int. Equ. \textbf{10}  (2007),
1131--1152.

\bibitem{Oh1} 
Y.-G. Oh,  \emph{Existence of semiclassical bound states of nonlinear Schrödinger equations with potentials of the class $(V)_a$}. Comm. Partial Differential Equations \textbf{13} (1988),  1499--1519.

\bibitem{Oh2}
Y.-G. Oh,  \emph{Correction to: ``Existence of semiclassical bound states of nonlinear Schrödinger equations with potentials of the class $(V)_a$''}. Comm. Partial Differential Equations \textbf{14} (1989),  833--834.

\bibitem{R}  P. H. Rabinowitz, \emph{On a class of nonlinear Schr\"{o}dinger equations.} Z. Angew. Math. Phys. \textbf{43} (1992),  270--291.



\bibitem{S} W. A. Strauss,  \emph{Existence of solitary waves in higher dimensions.} Comm. Math. Phys. \textbf{55} (1977),  149--162.


\end{thebibliography}
\end{document}